\documentclass[english]{ws-jktr}

\begin{document}

\title{The Classification of Spun Torus Knots}

\author{Blake Winter\\
\emph{State University of New York, Buffalo}\\
\emph{Email: bkwinter@buffalo.edu}}

\maketitle
\begin{abstract}
S. Satoh has defined a construction to obtain a ribbon torus knot
given a welded knot. This construction is known to be surjective.
We show that it is not injective. Using the invariant of the peripheral
structure, it is possible to provide a restriction on this failure
of injectivity. In particular we also provide an algebraic classification
of the construction when restricted to classical knots, where it is
equivalent to the torus spinning construction.
\end{abstract}
\section{Introduction}
Throughout this paper we will work only with oriented, 1-component classical
knots, oriented, 1-component welded knots, and oriented surface knots. Furthermore,
all isotopies are assumed to be orientation-preserving.\\
S. Satoh\cite{SS} has shown that there is an algorithm to produce an oriented ribbon
torus knot from any oriented welded knot diagram. We follow his notation and
designate this operation as $Tube$. Furthermore, this operation was
shown to be independent of the particular representative of a welded-equivalence
class of welded knots. This was proved by showing that any welded
Reidemeister move induces an isotopy on the corresponding ribbon torus
knots. Satoh also demonstrated that this operation is surjective,
in the following sense: for any isotopy class of ribbon torus knots,
there is some welded knot $K$ such that $Tube(K)$ lies in that isotopy
class.\\
It is natural to ask whether or not ribbon torus knots are classified
by welded knots under this $Tube$ operation; that is, if $Tube(K)$ and $Tube(L)$
are isotopic, must $K$ and $L$ be welded-equivalent? It will be shown that this is not the case,
by exhibiting a specific example of inequivalent welded knots which
are mapped to the same ribbon torus knot by $Tube$. We will consider knots which are not $(-)$ amphichiral,
and show that for such knots, $Tube$ fails to be injective.\\
We will then examine the peripheral structure for oriented, 1-component classical and welded knots,
and extend this invariant to surface knots. Using this it is possible to determine
that for classical knots, $Tube^{-1}(Tube(K))=\left\{K,-K^{*}\right\}$. This leads to an algebraic
classification of oriented spun classical knots.
\section{Welded Knots}
Kauffman\cite{Kauf} introduced virtual knots as a combinatoric generalization
of classical knot diagrams, in which a third crossing type is allowed, a virtual crossing. 
Welded knots\cite{FRR,SS,SK} are presented via
the same diagrams as virtual knots, but with one additional move, one of the 'forbidden'
moves from virtual knot theory (see Fig. 1). Two welded knot diagrams which may be transformed into each other by a sequence of such moves are called \emph{welded
equivalent} or \emph{w-}equivalent.\\
\begin{figure}
\begin{centering}
\includegraphics[scale=1]{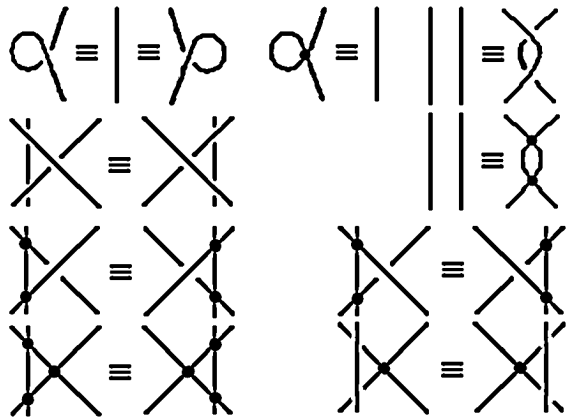}
\par\end{centering}

\caption{The welded Reidemeister moves. The move in the lower right is one of the
'forbidden' moves.}

\end{figure}
Any classical knot diagram may be interpreted as a welded
knot diagram. Our first theorem shows that if two classical knot diagrams
are welded equivalent, then they are equivalent as classical knots
as well. In order prove this, we first recall the definition of the knot group and
longitude for classical and welded knots. For both classical and welded knots,
the knot group may be found in a Wirtinger presentation from the diagram. Each
arc in the diagram is a generator, and each crossing yields a relation of the form
$x=z^{-1}yz$, which may be abbreviated as $x=y^{z}$.
\begin{figure}
\begin{centering}
\includegraphics[scale=1]{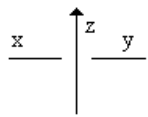}
\par\end{centering}

\caption{Each crossing in a knot diagram yields a relation of the form $x=y^{z}$. Reversing
the orientation of the overcrossing would change the relation to $y=x^{z}$.}

\end{figure}

The longitude is defined as an element of the knot group corresponding
to circling the knot exactly once in the direction of the
orientation of the knot without algebraic linking; it is the canonical generator of the non-linking part of the 
peripheral group\cite{DJ, GPV}. Thus, it is an ordered list of the arcs one crosses under,
where the arc's generator appears at right handed crossings and its inverse
appears at left handed crossings, multiplied at the end by
$m^{-k}$ where $k$ is the sum of the signs of the crossings, and 
$m$ is the meridian, so that the linking number is $0$.\\
This latter, combinatorial definition of the longitude extends naturally to welded knot
diagrams. The longitude of a welded knot is defined to be the element of the knot group obtained
by multiplying the generators of the arcs which one passes under in classical crossings,
multiplied at the end by $m^{-k}$ where $k$ is the sum of the signs
of the crossings, and $m$ is the meridian. Welded crossings do not contribute to the longitude.
We see that this definition is in fact invariant by considering their behavior
under the welded Reidemeister moves. For example, we may check that under the 'forbidden move,' which moves
an overcrossing over a welded crossing, the longitude does not change: the order in which the overcrossing will appear in the list remains unchanged.\\
Since the combinatorial definitions of the longitude are identical for both classical and welded knot diagrams, the longitude 
of a classical knot is the same whether computed using the welded definition
or the classical definition. For the order
in which we encounter overcrossings does not depend upon whether the knot diagram
is being considered as a classical diagram or a welded diagram.\\
It is known\cite{Wal,Kauf,DJ,CG} that the \emph{group system} (the knot group, the meridians, and their corresponding longitudes) classifies (oriented, 1-component) classical knots.\\
\begin{theorem} If $K$ and $L$ are classical (oriented, 1-component) knots whose diagrams
are welded equivalent, then they are isotopic.\end{theorem}
Proof (\cite{GPV,Kauf}): The group system is preserved
under welded Reidemeister moves.
Therefore $K$ and $L$ must have the same classifying invariant,
and be classically equivalent.$\square$
\section{Satoh's Tube Map}
For full proofs and exposition of the following, we refer the 
reader to Satoh's development\cite{SS}. We will review here only the essential points of
the construction.\\
Satoh defined the operation $Tube$ as follows. Given a welded knot
diagram $K$, we draw a broken surface diagram\cite{CS} by placing a thin
tube wherever we see an edge in the welded knot diagram, orienting
the surface as shown in Fig. 2. At welded crossings, the tubes pass
over/under each other, and at classical crossings, they knot together
as in Fig. 2. We orient the surface so that the normal vector in our
broken surface diagram points outward.\\
\begin{figure}
\begin{centering}
\includegraphics[scale=1]{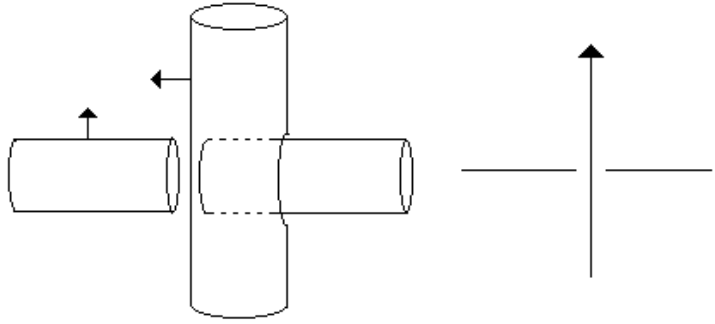}
\par\end{centering}

\caption{Satoh's definition of the Tube map.}

\end{figure}
The surface corresponding to the resulting broken surface diagram
is defined to be $Tube(K)$. We have from Satoh's work the following
important results.
\begin{theorem} For any two welded knots $K,\, K'$, if $K\cong K'$,
$Tube(K)\cong Tube(K')$.\end{theorem}
Satoh proves this by showing that any welded Reidemeister move can
be mirrored on the broken surface diagram using Roseman moves.
\begin{theorem} For any ribbon torus knot $R$, there is a welded
knot $L$ such that $Tube(L)\cong R$.\end{theorem}
In the following the $*$ operation refers to taking the mirror image
of the surface knot.
\begin{theorem} The following equivalences hold for any welded knot
$K$: $Tube(K)\cong -Tube(K)^{*}, -Tube(K)\cong Tube(-K)$.\end{theorem}
From this it follows that $-Tube(-K)\cong Tube(K).$
We recall here the operation on classical knots which Satoh denotes by
$Spun(K)$. Take a classical knot $K$ and place it in a half-hyperplane copy of
$\mathbb{R}^{3+}$ within $\mathbb{R}^4$. Now rotate the half-hyperplane about its
face. $K$ will trace out a surface with the diffeomorphism class of
a torus. We refer to this torus as $Spun(K)$.
\begin{theorem}For oriented classical knots $K$, $Tube(K)\cong Spun(K)$ for at least one orientation
of $Spun(K)$, and exactly one if $Spun(K)$ is not reversible.
We denote this orientation of $Spun(K)$ by $OSpun(K)$; that is,
$OSpun(K)$ is the orientation of $Spun(K)$ which makes $Tube(K)\cong OSpun(K)$ true.\end{theorem}
Proof: We refer the reader to Satoh's construction of an explicit
isotopy of the unoriented $Spun(K)$ to the (unoriented) $Tube(K)$.
By requiring that $OSpun(K)\cong Tube(K)$ as oriented surfaces, we
induce an orientation on $OSpun(K)$ via this isotopy, which is clearly
well defined, unless $Spun(K)$ is reversible, in which case $OSpun(K)$
is reversible as well, and hence it is well defined in this case as
well. Note that if $Spun(K)$ is not reversible, then neither is $OSpun(K)$
nor $Tube(K)$, for if they were Satoh's construction could be applied
to reverse the orientation of $Spun(K)$.$\square$\\
\begin{theorem} $Tube$ preserves the knot group and quandle.\end{theorem}
This follows from a straightforward computation (see the diagram).
\begin{figure}
\begin{centering}
\includegraphics[scale=1]{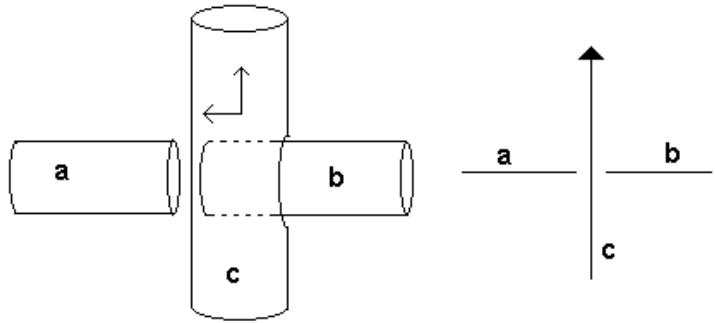}
\par\end{centering}
\caption{Computing the relation from the tube diagram, we obtain
$a=b^{c}$, just as we obtain from the welded diagram. It follows that
the quandle presentations we obtain are identical.}
\end{figure}
We require an additional theorem, which slightly generalizes a theorem
pointed out to the author by Dennis Roseman. Recall that vertical
reflection of a welded knot diagram is performed by reflecting the planar graph across a vertical plane\cite{HK}. For
classical knots this reduces to the usual reflection. We denote the vertical
reflection of $K$ by $K^{\uparrow}$.
\begin{theorem}For welded knots $K$, $Tube(K)^{*}\cong Tube(-K^{\uparrow})$.\end{theorem}
Proof: Draw the diagram of $Tube(K)$ and then draw the mirror, $Tube(K)^{*}$.
Now look at the welded knot $K'$ which naturally yields $-Tube(K)^{*}$. 
This is precisely $-K^{\uparrow}$, so $-Tube(K)^{*}\cong Tube(-K^{\uparrow})$.$\square$
\begin{figure}
\begin{centering}
\includegraphics[scale=0.4]{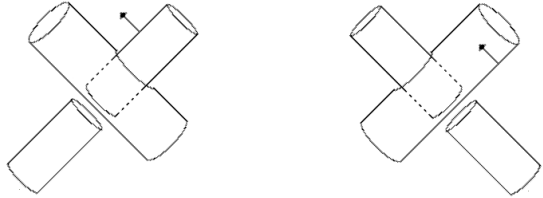}
\par\end{centering}

\caption{A crossing of tubes and its mirror image.}

\end{figure}
\subsection{Tube is Not Injective}
We are now ready to prove our main theorem about Satoh's construction.
In private correspondence Satoh has indicated that he became aware
of this theorem, or the possibility of it, at some point after \cite{SS} was published.
\begin{theorem} There is a ribbon torus knot $R$ with the property
that $Tube^{-1}(R)$ contains inequivalent welded knots.\end{theorem}
Proof: Let $K$ be a welded knot, with $K\ncong -K^{\uparrow}$.
Such knots exist; take for instance the right handed trefoil which is chiral as a welded knot, since classical
knot theory embeds faithfully in welded knot theory. But we have shown that
$Tube(K)\cong -Tube(K)^{*}\cong Tube(-K^{\uparrow})$, so $Tube^{-1}(Tube(K))$ contains the
inequivalent welded knots $K$ and $K^{\uparrow}$.$\square$\\
Satoh has also defined \emph{arc diagrams}\cite{SS}, which are welded knots with endpoints, in order to describe ribbon 2-knots (which have the diffeomorphism type of a sphere). However, it is not clear whether this theorem
extends to this case or not. In particular, the theorem must extend if there are
arc diagrams which are not $(-)$ amphichiral.\\
It remains an open question to determine the extent
to which $Tube$ fails to be injective. For example, is it possible
to place an upper bound on the cardinality of $Tube^{-1}(R)$, or to
describe this set precisely? By considering
the peripheral structure we can place a partial bound on this set, and 
also show that the peripheral structure is a classifying invariant
on the subset of ribbon torus knots consisting of oriented spun tori.
To prove this we will need to carefully examine the peripheral group\cite{Kauf,DJ,CG} 
classical knots.\\
\section{The Longitude Group Invariant}
Recall that a surface knot is an embedding of a surface, of arbitrary genus,
into $S^{4}$. We define the \emph{longitude group} invariant generally, for any surface knot regardless
of genus, in a way which generalizes the longitude invariant of a
classical knot.
\\Define the \emph{linking number} of a loop with a surface knot
as follows. Given a Wirtinger presentation of the knot group $G$ from the broken surface
diagram of a surface knot there is a homomorphism
$G\rightarrow \mathbb{Z}$ defined as the sum of the exponents in a word defining
the element. This is well defined, as any changes to the word involve replacing
a generator with a conjugate of a generator, or the reverse, and therefore do not
change the sum of the exponents. Similarly, under Roseman moves on the diagram,
the definition remains unchanged.\\
Let $R$ be a surface knot of genus $n$ embedded
in a manifold $M$. Let $N(R)$ be a regular neighborhood of $R$, whose closure is contained in another regular neighborhood of $R$, and let $X=M-N(R)$. Let $B=\partial X$. Observe that $B=\partial\overline{N(X)}$, as well. There is a natural embedding $i:B\hookrightarrow M$, which is the inclusion. The embedding induces
a homomorphism $i_{*}:\pi_{1}(B)\rightarrow\pi_{1}(X)$. We refer
to the image of this homomorphism as a \emph{peripheral group} of
$R$, $P(R)$. The longitude group is now defined to be those elements in the
image which do not link with the knot. From the definition, the linking
number is additive under composition, from which it follows that
this set is in fact a group. We define this to be the \emph{longitude group}, and denote it
by $LG(R)$.\\
Now any element of the knot group for an oriented surface represents some homology
class of curves. The first homology will be isomorphic to the integers, since the
abelianization of any group with a Wirtinger presentation, such as the knot group,
is isomorphic to the integers. The \emph{linking number} of a curve representing
an element in the knot group is defined to be the image in the integers of its
homology class under the natural isomorphism.\\
We define an element of $P(R)$
with linking number $1$ to be a meridian (of the knot group) if it bounds a hyperdisk inside $B$ about a point on the surface knot and has linking number 1. The meridian element of the peripheral group is denoted by $m(R)$. Uniqueness for
the meridian follows from the fact that any regular neighborhood of an oriented knotted surface is
a trivial disk bundle.\\
The triple $(P(R),LG(R),m(R))$ will be called the peripheral structure.
Note also that this definition may be immediately generalized to other dimensions. 
However except in dimension $1$ the longitude group may be acyclic, and even if it is cyclic
there need not be a canonical generator.\\
Observe that one may choose a different meridian, which defines a different longitude group
conjugate to the first one.
We will use the term \emph{peripheral strucure} to refer
to this collection and its conjugacy relations. Recall that the quandle
of a surface knot may be defined in terms of nooses which link
with the surface, as Joyce defined them for one dimensional knots\cite{DJ}. By the topological definition of all these
constructions,
\begin{theorem} If $R,R'$ are surface knots which are ambiently
isotopic, then there is an isomorphism of their quandles, which induces
an isomorphism of their peripheral structure.\end{theorem}
We can also extend a theorem of Joyce about quandles of classical knots and peripheral groups. Recall
that the knot group has a natural right action on the quandle
given by composing the quandle noose with the group loop.
\begin{theorem} For a surface knot $R$ with group $G$ and quandle $Q$, let $q\in Q$, and let
$G_{q}=\lbrace g\in G \vert qg = q\rbrace$. Then $G_{q}$ is a peripheral group of $R$.
Conversely any peripheral group stabilizes some quandle element.\end{theorem}
Proof: See \cite{DJ}, Thm. 16.1; the proof generalizes identically. We sketch the proof here. Given
a peripheral group, choose a quandle element $q$ such that the meridian
is representable by a path which follows the quandle noose and then comes
back upon it. Then the meridian stabilizes $q$.
Now given an element $g$ of the longitude group, we may perform a homotopy so that this element lies
in a regular neighborhood of a surface. Acting on $q$ with this element, we may now
pull the head of $q$ along a path on the surface parallel to $g$, which will end back
at $q$ itself.\\
The other direction is similar. Given $q$ we choose the peripheral group
whose meridian is parallel to $q$. Then that peripheral group stabilizes $q$. No other elements
stabilize it, however. For if some element $g$ stabilizes it, then the homotopy from $qg$ to $q$
implies that we may push $g$ via a homotopy to make it into a path lying
completely within some regular neighborhood of the surface conjugated by $q$. 
It must therefore be in the peripheral group with meridian $q$. See \cite{DJ} for details.$\square$\\
\subsection{Computing Longitude Groups for Surface Knots}
Consider an arbitrary surface knot $R$. Let us take a broken surface 
diagram (although the method will generalize to any desired presentation).
Observe that the Wirtinger generators, those which loop through only
one surface in the broken surface diagram, each have linking number $1$.\\
Consider the boundary of a regular neighborhood, $B$. The regular neighborhood
will be a trivial disk bundle $D^{2} \times R$, as $R$ has codimension two and all orientable
surfaces in $S^{4}$ have normal Euler number $0$. Therefore the boundary of this disk bundle
will be $B\cong R \times S^{1}$. Let us take a particular generating set
for $\pi _{1}(B)\cong \pi _{1}(R)\times \mathbb{Z}$ consisting of $\lbrace a_{1},...,a_{n}\rbrace$ the
generators of the $R$ component, and $m$ the generator of the $S^{1}$
component. We choose these such that $i_{*}(m)$ is the meridian, with
linking number 1, and such that to find $i_{*}(a_{j})$ we take a
path on the surface of the broken surface diagram and perturb it
off the surface. It follows then that to write down $i_{*}(a_{j})$ 
we need only record the signed overcrossings which we encounter
as we loop around the broken surface diagram.\\
Now we wish to change to a collection of generators which will
make the longitude group and the meridian easily recognizable. To do so, let $k_{j}$ be the sum of the
signs of the overcrossings encountered along $a_{j}$; that is,
it is the linking number of $a_{j}$. Let us construct a new
generating set for $\pi _{1}(B)$ consisting of $\lbrace m, b_{1},...,b_{n}\rbrace$,
where $b_{j}=a_{j}m^{-k_{j}}$. It is straightforward to verify
that this really is a generating set. Furthermore, all the images of the generators have
linking number 0 except for $i_{*}(m)$ which has linking number 1.
But this is a generating set of the peripheral group. Therefore, $LG(R)$ is generated by
$\lbrace i_{*}(b_{1}),...,i_{*}(b_{n})\rbrace$, since these have linking number $0$.\\
Now notice also that $P(R)$ is the image of a homomorphism
from $\pi _{1}(R)\times \mathbb{Z}$, and that $m(R)$ is the image of
the generator of $\mathbb{Z}$, which commutes with anything in $\pi _{1}(R)$ part of
the product. Therefore,
\begin{theorem} For any surface knot $R$ with the diffeomorphism
type of $N$, $LG(R)$ is a quotient of $\pi_{1}(N)$,
and $m(K)$ commutes with all elements of the peripheral group.\end{theorem}
Note that this implies that the longitude group of any sphere-knot
is trivial.\\
We now restrict our attention to the specific case of tori. For a
torus, $B\cong S^{1}\times S^{1}\times S^{1}.$ Following the above
computation we see that the longitude group
is isomorphic to a quotient of $\mathbb{Z}\oplus \mathbb{Z}$.\\
\begin{figure}
\begin{centering}
\includegraphics[scale=0.4]{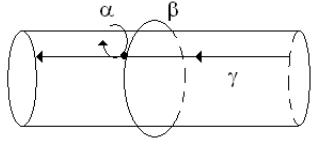}
\par\end{centering}

\caption{The three generators of the peripheral group for a torus.}

\end{figure}
\begin{theorem} If $R$ is a surface knot with the diffeomorphism
class of a torus, then $LG(R)$ is a quotient of $\mathbb{Z}\oplus\mathbb{Z}$.\end{theorem}
For welded knots, define the peripheral group to be the subgroup
of the knot group generated by the longitude and meridian. Observe that this
combinatorial definition agrees with the topological definition for classical knots.
Now consider $R\cong Tube(K)$ where $K$ is a welded knot. Then we have the following
important result about Satoh's construction.
\begin{theorem} The peripheral group, meridian, and longitude
group are preserved by $Tube$.\end{theorem}
Proof: Let the generators of $\pi _{1}(B)$ be $\alpha ,\beta ,\gamma$.
Map $\alpha$ to the meridian of the peripheral group. We may choose
this to coincide with the meridian of the welded knot, which proves the second claim.
We send $\beta$ to the path which loops around the meridian of the torus
(not to be confused with the meridian of the peripheral group),
which is a trivial loop in the knot group and therefore
already has linking number $0$. Finally we send $\gamma$ to the path
traveling around the equator of the torus. $i_{*}(\gamma)$ is therefore
given by recording the signed overcrossings it encounters. We take a new
generating element $\delta = \gamma \alpha^{-k}$ where $k$ is the sum
of these signs. As before, $\lbrace \alpha ,\beta ,\delta \rbrace$
form a generating set for $\pi_{1}(B)$, and furthermore $LG(R)$ is
generated by $i_{*}(\delta)$, since any element in $LG(R)$ must have
a reduced word consisting only of $i_{*}(\delta)$ and $i_{*}(\beta)$, the latter
being trivial. But $i_{*}(\delta)$ is also the longitude
of $K$, or possibly its inverse. The cyclic subgroup generated by $i_{*}(\delta)$
is the longitude group of $K$ in either case. Therefore the longitude group
is preserved, and since the meridian is also preserved, the peripheral group is as well.$\square$\\
However, it is not possible to choose a preferred generator
of the longitude group for ribbon tori, as one may for classical knots
(for if it were possible, then $Tube$ would be injective
when restricted to classical knots). In particular, for classical knots
the orientation induces a choice of canonical generator for the longitude group,
whereas the orientation of a surface does not. However, there are only
two choices for generators of the longitude group for the surface, $i_{*}(\delta)$ and $i_{*}(\delta)^{-1}$.
Our computation shows that one of these is the longitude of the corresponding welded knot,
while the other is its inverse. Note also that
as a consequence of the above theorem and the fact that 
the peripheral group of any torus (surface) knot is abelian (being
the image of a homomorphism of $\pi _{1}(S^{1}\times S^{1}\times S^{1})$) that the peripheral group of a welded
knot is abelian.\\
Define the peripheral
structure of a welded knot to include all the
peripheral groups and their conjugation relations;
from these definitions it follows that
\begin{theorem} $Tube$ preserves the conjugacy relations among the peripheral groups, and hence
the peripheral structure.\end{theorem}
As an immediate consequence from our computation we have
\begin{theorem} If a torus embedding $R$ is ribbon, then $LG(R)$
is cyclic.\end{theorem}
However, this condition is not sufficient. Consider the so-called "1-turned trefoil torus," which is studied
in \cite{Boy}.
This is constructed from a trefoil by rotating the classical knot around a $S^{2}\subset S^{4}$ surface,
while also turning it $2\pi$ around another $S^{2}\subset S^{4}$. The resulting longitude group
is isomorphic to that of the trefoil and hence cyclic. To compute this
we choose one generator to be the one which is parallel to the trefoil itself,
and the other to be one which spins and rotates with some point on the trefoil. The latter
loop is contractible, whereas the former generates the longitude group of the 1-turned trefoil torus\cite{Tan}. However this torus embedding fails to be ribbon \cite{Boy}.
Additionally not all nontrivial tori will have nontrivial longitude groups; indeed the connected sum of a trivial torus embedded in $S^{4}$ with any knotted sphere will have trivial longitude group. For the generators of the longitude group
will be homotopic to paths parallel to the meridian and longitude of the trivial torus, both of which are contractible.
To summarize our above results, then:
\begin{theorem} $Tube$ preserves the knot group, the meridian, and the longitude up to inverse, as well as the conjugacy
relations between different peripheral groups.\label{TubePreserve}\end{theorem}
\section{A Classification of Oriented Spun Tori}
Throughout this section we will work only with oriented 1-component classical knots.
We will therefore denote the mirror image of a knot $K$ by $K^{*}$, since there is
no ambiguity.\\
The triple $(\pi _{1}, P, m)$ (hereinafter referred to as
the peripheral structure) classifies classical oriented knots up to mirror reverses\cite{Wal,Neu,Man}.
\begin{theorem} If $K$ and $K'$ are oriented one-component classical knots with isomorphic
peripheral structures, then either $K \cong K'$ or $K \cong -K'^{*}$.\end{theorem}
From this theorem, together with Thm. \ref{TubePreserve}, we have
\begin{theorem} For classical oriented knots $K,K'$, $Tube(K)\cong Tube(K')$ (or equivalently $OSpun(K)\cong OSpun(K')$) iff
$K\cong K'$ or $K\cong -K'^{*}$.\end{theorem}
Proof: If $Tube(K)\cong Tube(K')$ then $K$ and $K'$ have isomorphic peripheral structures,
and therefore $K\cong K'$ or $K\cong -K'^{*}$. On the other hand if $K\cong K'$
then $Tube(K)\cong Tube(K')$ automatically. If $K\cong -K'^{*}$, then
$Tube(K')\cong -Tube(K')^{*}\cong Tube(-K'^{*})\cong Tube(K)$.$\square$\\
As a consequence, we have an algebraic version of this theorem.
\begin{theorem} Oriented spun tori are classified by their peripheral structures.\end{theorem}
It is an open question to determine whether the peripheral structure suffices
to classify ribbon torus knots generally. It is known that it does not classify torus
embeddings up to isotopy generally\cite{Tan}. A related problem is to determine
whether the knot group, meridian, and longitude of an oriented welded knot
suffice to determine the welded knot. A positive answer to the latter problem
would imply a positive answer to the former question, using the same method
as used above to classify the oriented spun tori. Conversely, if the peripheral
structure classifies ribbon tori, this would imply that it classifies oriented welded knots
up to reversed vertical reflection.\\
\paragraph{Acknowledgements}
The author would like to thank Sam Nelson, Louis Kauffman, Dennis
Roseman, and William Menasco for their discussions and suggestions,
Xiao-Song Lin for his discussions and for alerting the author
to the problem, and the referee for many helpful comments on
the paper itself.

\end{document}